\pdfoutput=1

\documentclass[preprint, review, 12pt, 3p]{elsarticle}

\usepackage{amssymb, csquotes, listings, adjustbox, longtable, tabu, amsmath, mathtools, enumitem}

\usepackage{hyperref}

\usepackage{floatpag}
\floatpagestyle{empty}

\journal{}

\begin{document}

\begin{frontmatter}

\vspace*{1.2in}

\title{Geometry of Arrangements that Determine Shapes}

\author[label1]{Alexandros Haridis}
\address[label1]{Department of Architecture, Massachusetts Institute of Technology}
\ead{charidis@mit.edu}

\begin{abstract}
Shape grammars compute over shapes which are defined in the universe $U^*$. Shapes in the universe $U^*$ are analogous to line drawings that can be physically realized in the plane. In this universe, a shape is embedded or contained in an \emph{arrangement} of ``points" and ``lines" in the plane called, respectively, \emph{registration marks} and \emph{construction lines}. In this expository article, arrangements that contain shapes are studied as incidence structures and the finite geometries they give rise to are characterized. In particular, arrangements that contain shapes are distinguished into those that give rise to finite near-linear and linear spaces, and those that do not give rise to any proper form of ``geometry." Arrangements that constitute finite geometries (near-linear and linear spaces) give an alternative characterization of determinate rules in shape grammars. This paper contributes to the body of work related to the \emph{mathematics of shapes} in the area of shape grammar theory.
\end{abstract}


\end{frontmatter}



\vspace{2.5in}

\section{The universe of shapes $U^*$}
\label{sec1}

A \emph{shape} is given by a finite set of line segments in the two-dimensional Euclidean plane:

\vspace{12pt}
\centerline{$\{\;l_1, ..., l_n\;\}$}
\vspace{12pt}

\noindent such that there are no line segments $l_1, ..., l_k$, $1 < k \leq n$, for which $l_1\;\mathtt{co}\;l_2$, $l_2\;\mathtt{co}\;l_3$, ..., $l_{k-1}\;\mathtt{co}\;l_k$. The relation $\mathtt{co}$ \cite[p.~127]{StinyPhD1975}, identifies sets of line segments that are adjacent co-linear or overlapping co-linear; see Figure 1. Thus, a shape consists of a finite number $n \geq 0$ of line segments in the plane, so that no two line segments are adjacent co-linear or overlapping co-linear. A shape that contains no line segments ($n = 0$) is denoted by $S_{\emptyset}$ and is called \emph{empty shape}.

\begin{figure}[b!]
\centering
\includegraphics{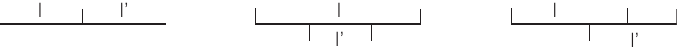}
\caption{Adjacent co-linear and overlapping co-linear line segments are not permitted to form shapes.}
\label{Figure1}
\end{figure}

Given two shapes $S_1$ and $S_2$, their union $S_1 \bigcup S_2$ is given by the union of their respective sets of line segments. For a finite index set $\mathcal{J}$ and a set of shapes $\mathcal{S} = \{\;S_i\;|\;i \in \mathcal{J}\;\}$, the union of the shapes in $\mathcal{S}$ is given by $\bigcup_{i \in \mathcal{J}}S_i$. 

The union of two or more shapes may not necessarily be a shape. The resulting set may contain adjacent co-linear or overlapping co-linear line segments (Figure 2). 

To obtain a shape, the set resulting from the union must be reduced to a set of maximal line segments: the smallest possible set of line segments that is pictorially equivalent to the union, and does not contain line segments like those in Figure 1. In particular, the set $(S_1 \bigcup S_2)^R$ is the shape formed from the reduced union of $S_1$ and $S_2$. The definition of the reduced union of two shapes is not constructive. The reduced union is instead obtained algorithmically, using reduction rules for line segments \cite{StinyPhD1975}.

\begin{figure}[t!]
\centering
\includegraphics{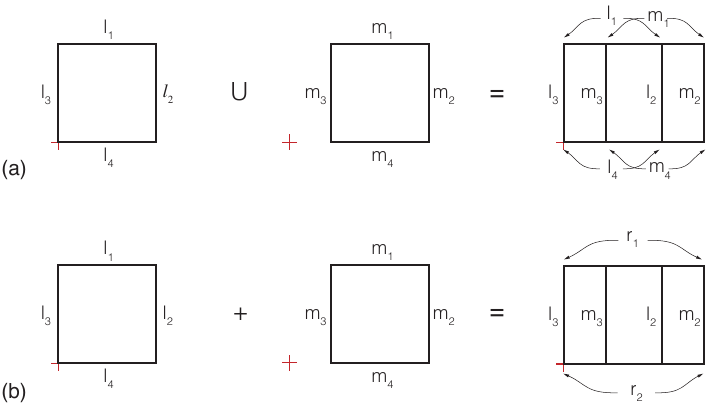}
\caption{(a) The union of two shapes may not be a shape, since the result may contain adjacent co-linear or overlapping co-linear line segments. (b) The sum of two shapes is always a shape. The red cross-hair indicates the origin of the plane.}
\label{Figure2}
\end{figure}

The sum (+) of two shapes $S_1$ and $S_2$, denoted by $S_1 + S_2$, is given by:

\vspace{12pt}
\centerline{$S_1 + S_2 \coloneqq (S_1 \bigcup S_2)^R$}
\vspace{12pt}

\noindent When the operation of sum is combined with transformations or sequences of transformations, the universe of all shapes is obtained. 

In particular, if $\mathcal{S}$ is a finite, non-empty set of shapes, then the set $\mathcal{S}^+$ is defined as the smallest set containing the shapes in $\mathcal{S}$ and closed under sum and the transformations \emph{translation}, \emph{rotation}, \emph{reflection about an axis}, and \emph{uniform scale}. A shape is in $\mathcal{S}^+$ when it is in $\mathcal{S}$ or it can be formed by applying sum or the transformations a finite number of times to a shape or shapes in $\mathcal{S}$. For a set of shapes $\mathcal{S}$, the set $\mathcal{S}^*$ is the union of $\mathcal{S}^+$ and a set containing just the empty shape, i.e. $\mathcal{S}^* = \mathcal{S}^+ \bigcup \{\;S_{\emptyset}\;\}$. A shape is in $\mathcal{S}^*$ if it is in $\mathcal{S}^+$ or it is the empty shape.

Now, let $U$ be a set that contains a shape consisting of a single line segment $\{((0, 0), (1, 0))\}$. Then, $U^*$ = $\{\;\{((0, 0), (1, 0))\}\;\}^*$ is the set of all possible shapes and it is called \emph{universe of shapes} \cite{StinyPhD1975}. The universe contains:

\vspace{9pt}

\begin{quotation}
\noindent All possible finite sets of line segments in the plane, such that no two line segments in each set are adjacent co-linear or overlapping co-linear.
\end{quotation}

\noindent Since the transformations are defined over real parameters, the universe is uncountable.

\section{A shape is contained in a point-line arrangement}
\label{sec2}

Each line segment $l = (p_1, p_2)$ of a shape has an equation $eq(l)$ which is determined by its two end-points, $p_1$ and $p_2$. The equation of a line segment represents an ``infinitely extending line" that contains the segment (and possibly other co-linear line segments that are separated by ``gaps").

Let \emph{L} be a set of infinitely extending lines containing the segments of a shape $S$:

\vspace{12pt}
\centerline{$L = \{\;eq(l) \;|\; eq(l)\;\textrm{is the equation of a line containing}\; l \in S\;\}$}
\vspace{12pt}

\noindent The set \emph{L} is called \emph{construction lines} of $S$. For convenience, the equation $eq(l)$ of a line is denoted by $\mathcal{L}$.

A point $q$ is a \emph{point of intersection} or \emph{registration mark} if and only if there are construction lines $\mathcal{L}_1$ and $\mathcal{L}_2$ in $L$, such that $\mathcal{L}_1 \neq \mathcal{L}_2$, and $\mathcal{L}_1$, $\mathcal{L}_2$ intersect or meet at $q$. That is to say, the point $q$ is a solution to the equation of both lines. Two distinct construction lines do not meet if and only if they are parallel.

The set $P$, given by:

\vspace{12pt}
\centerline{$P = \{\;q \;|\; q \; \textrm{is a point of intersection of}\; S\;\}$}
\vspace{12pt}

\noindent is called the \emph{registration marks} of the shape $S$. The pair $\mathcal{A}_S = (P, L)$ is called the \emph{arrangement} associated with the shape $S$. The definitions of the sets of construction lines and registration marks are constructive. Using the two definitions, it is straightforward to compute the arrangement of any given shape in the universe.

It is important to have a few simple rules that allow one to decide when a set of points and a set of lines in the Euclidean plane form an arrangement for a shape.

A pair $(P, L)$ of a set $P$ of points and a set $L$ of straight lines in the Euclidean plane forms an arrangement for a shape, if and only if:

\begin{enumerate}
    \setlength{\parskip}{0pt}
    \setlength{\itemsep}{0pt plus 1pt}
    \item[(1)] any point in $P$ is a point of intersection of at least two lines in $L$, and 
    \item[(2)] any two nonparallel lines in $L$ have a unique point of intersection in $P$.
\end{enumerate}

The two rules or axioms (1) and (2), do not assume that there are any lines at all. In this case, there cannot be points, too, and we get a special object called the \emph{empty arrangement}, denoted by $\mathcal{A}_{\emptyset}$---the arrangement associated with the empty shape. If there are lines, the rules do not imply that there are points, unless two or more lines meet. On the other hand, if there are points, then there must be lines.

Figure 3 shows some shapes and their underlying arrangements. The dashed lines represent the lines in the set $L$ and the dots the points in $P$. Figure 4 shows sets of points and/or lines which are not arrangements for shapes (they fail to satisfy either (1), (2), or both). Many configurations of points and lines familiar from classical geometry, are not actually arrangements for shapes; for example, those in Figure 5.

\begin{figure}[t!]
\centering
\includegraphics{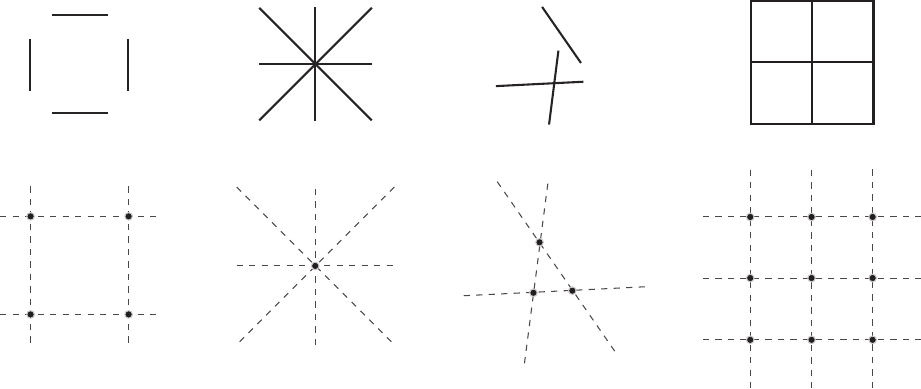}
\caption{Some shapes and their underlying arrangements.}
\label{Figure3}
\end{figure}

\begin{figure}[b!]
\centering
\includegraphics{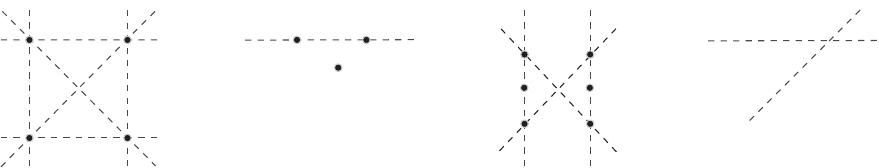}
\caption{Some sets of points and/or lines in the plane which are not arrangements for shapes.}
\label{Figure4}
\end{figure}

\begin{figure}[h!]
\centering
\includegraphics{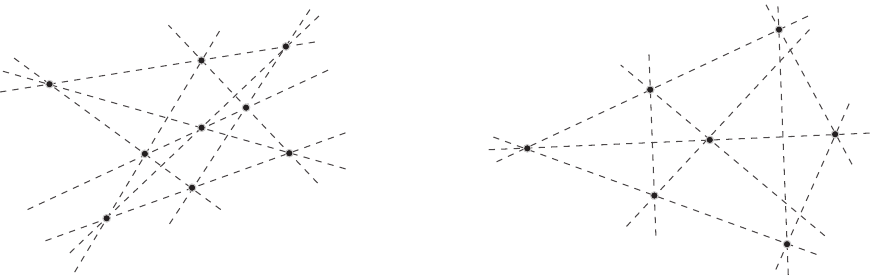}
\caption{Configurations of points and lines like Pappus's (left) and Desargues's (right) are not arrangements that determine shapes.}
\label{Figure5}
\end{figure}

\section{Operations between arrangements}

Apart from sum, shapes in the universe can combine in other operations, too, such as difference and product. Just like with sum, difference and product are also defined algorithmically \cite{StinyPhD1975}. When two shapes combine in an operation (sum, difference, or product) the result of the operation is always a shape. The arrangements of the shapes participating in an operation, are expected to combine in an analogous operation and to produce an arrangement as a result of that operation (i.e. to preserve the object type).

One way to define operations between arrangements is to focus on the set-theoretical definition of an arrangement and define the operations as operations between sets. Such an approach, however, quickly presents various complications.

In particular, given two arrangements $\mathcal{A}_{S_1} = (P_1, L_1)$ and $\mathcal{A}_{S_2} = (P_2, L_2)$, their union $\mathcal{A}_{S_1} \bigcup \mathcal{A}_{S_2}$ is given by the union of their respective sets of registration marks and construction lines, in particular, $P_1 \bigcup P_2$ and $L_1 \bigcup L_2$. The pair $(P_1 \bigcup P_2, L_1 \bigcup L_2)$, however, does not necessarily form an arrangement for a shape.

\begin{figure}[t!]
\centering
\includegraphics{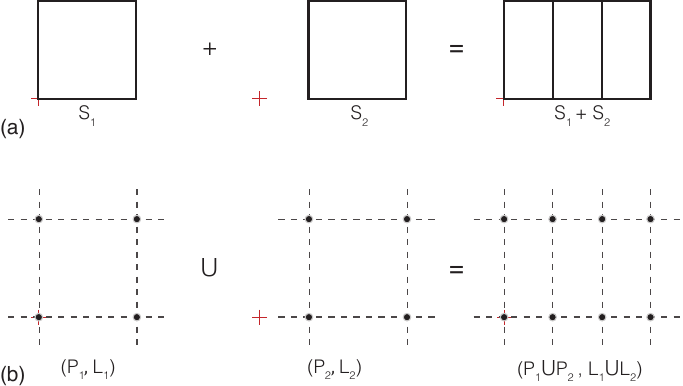}
\caption{}
\label{Figure6}
\end{figure}

In Figure 6, for example, the union of the two arrangements produces an arrangement---the resulting object satisfies the two rules (1) and (2), in section 2. For each sum in Figure 7, on the other hand, the union fails to produce an arrangement. Each union misses intersection points of construction lines (for economy of space, the lines in Figure 7 are not extended). The second union in Figure 7, for example, misses the two points colored in red in Figure 8.

\begin{figure}[b!]
\centering
\includegraphics{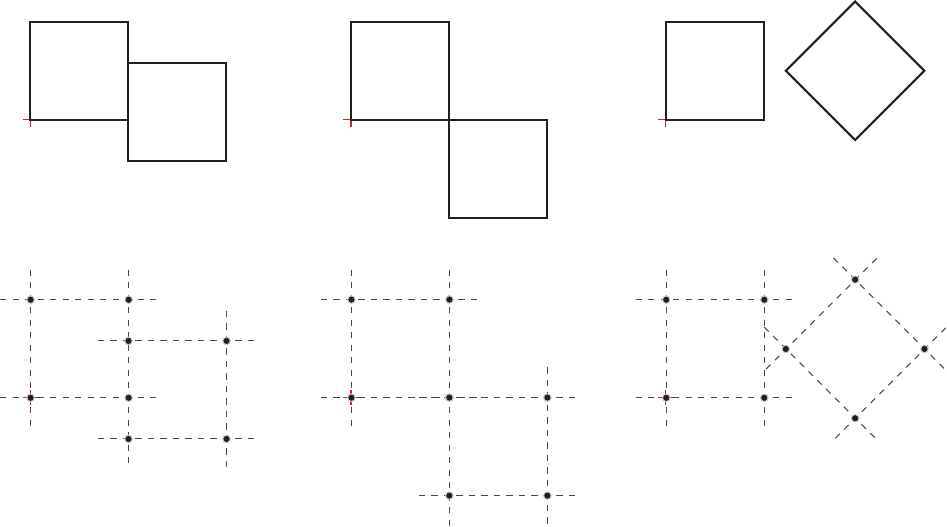}
\caption{}
\label{Figure7}
\end{figure}

\begin{figure}[t!]
\centering
\includegraphics{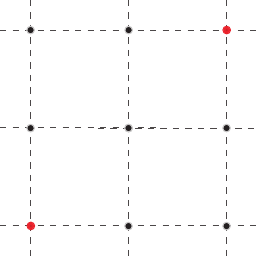}
\caption{}
\label{Figure8}
\end{figure}

Similar problems occur in other operations, too. The (set) difference of two arrangements, $\mathcal{A}_{S_1} \setminus \mathcal{A}_{S_2}$, is given by the set difference (relative complement) of their respective sets of registration marks and construction lines, in particular, $P_1 \setminus P_2$ and $L_1 \setminus L_2$. The pair $(P_1 \setminus P_2, L_1 \setminus L_2)$, however, does not necessarily form an arrangement for a shape.

\begin{figure}[b!]
\centering
\includegraphics{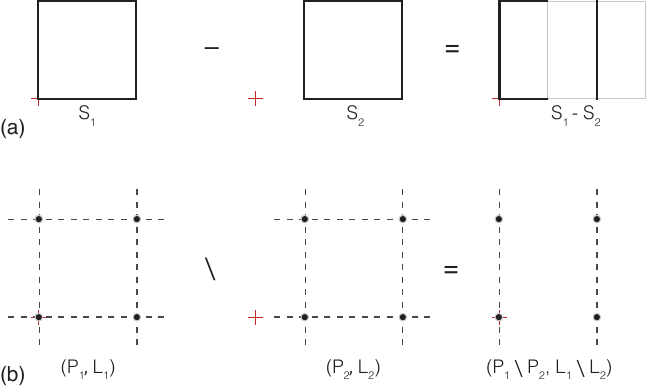}
\caption{}
\label{Figure9}
\end{figure}

In Figure 9, for example, the difference of the two arrangements fails to produce an arrangement. The resulting object does not satisfy rule (1). The same is true for the first two examples in Figure 10. The third example in the same figure satisfies rule (1), but not rule (2).

\begin{figure}[t!]
\centering
\includegraphics{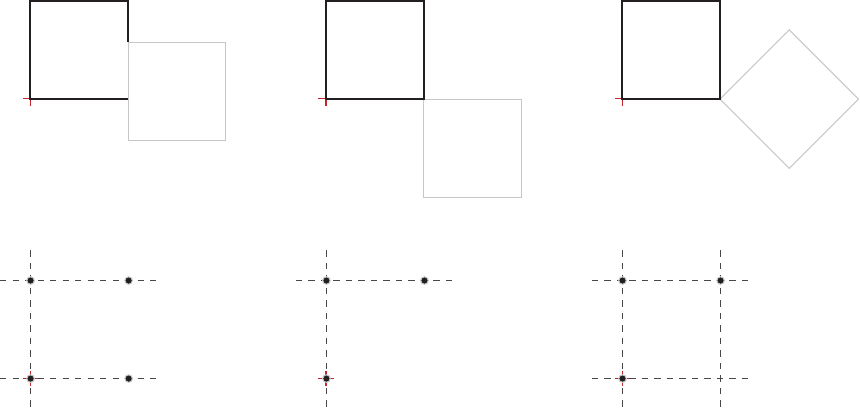}
\caption{}
\label{Figure10}
\end{figure}

The intersection of two arrangements, $\mathcal{A}_{S_1}$ and $\mathcal{A}_{S_2}$, is given by the intersection of their respective sets of registration marks and construction lines, in particular, $P_1 \bigcap P_2$ and $L_1 \bigcap L_2$. Again, the pair $(P_1 \bigcap P_2, L_1 \bigcap L_2)$ does not necessarily form an arrangement for a shape.

In Figure 11, for example, the intersection of the two arrangements produces the desired arrangement. On the other hand, the first two examples in Figure 12 are valid arrangements but they do not describe the desired shapes. The third example, in the same figure, is not an arrangement (it fails to satisfy rule (1)).  

\begin{figure}[b!]
\centering
\includegraphics{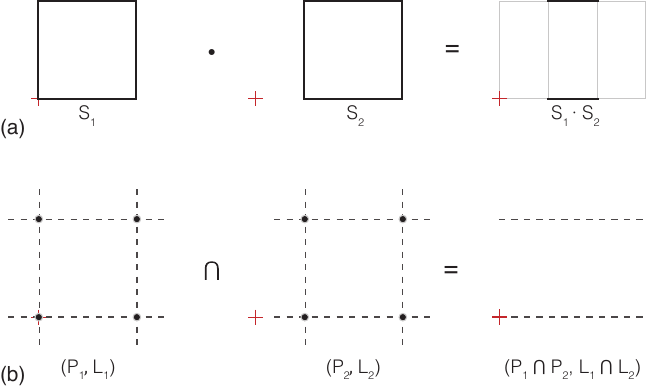}
\caption{}
\label{Figure11}
\end{figure}

\begin{figure}[t!]
\centering
\includegraphics{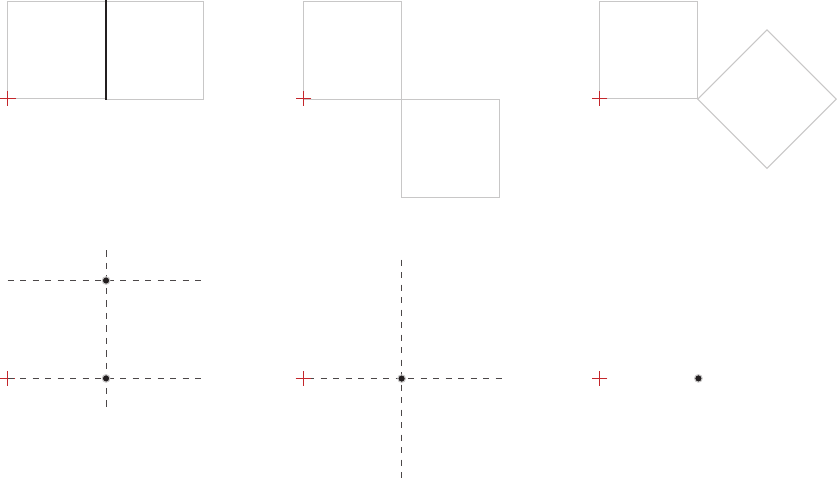}
\caption{}
\label{Figure12}
\end{figure}

It is possible to construct ``correcting rules" for each of the three operations, to fix the various problems case by case. These rules would essentially correct a resulting object obtained from an operation, whenever the object is not a valid arrangement. Such an approach, however, would require ad hoc decisions for how to handle every special case. 

Instead, a more natural way to define operations between arrangements is to define them \emph{relative} to the results obtained from the corresponding shape operations. 

In particular, given two shapes $S_1$ and $S_2$, the union of the arrangements $\mathcal{A}_{S_1}$ and $\mathcal{A}_{S_2}$, is defined as the arrangement of the shape $S_1 + S_2$. In notation,

\vspace{12pt}
\centerline{$\mathcal{A}_{S_1} \bigcup \mathcal{A}_{S_2} \coloneqq \mathcal{A}_{S_1 + S_2}$}
\vspace{12pt}

Analogous definitions are given for the difference and the intersection of two arrangements. In particular,

\vspace{12pt}
\centerline{$\mathcal{A}_{S_1} \setminus \mathcal{A}_{S_2} \coloneqq \mathcal{A}_{S_1 - S_2}$}

\vspace{12pt}
\centerline{$\mathcal{A}_{S_1} \bigcap \mathcal{A}_{S_2} \coloneqq \mathcal{A}_{S_1 \cdot S_2}$}
\vspace{12pt}

Given in this way, the definitions of the union, difference, and intersection of two arrangements are not constructive. However, the union, difference, and intersection of two arrangements can be obtained algorithmically, by determining the sets of construction lines and registration marks for, respectively, the shapes $S_1 + S_2$, $S_1 - S_2$, and $S_1 \cdot S_2$.

\section{Arrangements for shapes as incidence structures}

Incidence relations between registration marks and construction lines can be captured by an abstract mathematical object known as ``incidence structure". Incidence structures are at the core of the study of finite geometries, and their application extends to computer science, experiment design in statistics, and other areas (\cite{BattenCombinatoricsFiniteGeom}, \cite{Grunbaum2009}, \cite{Dembowski1997}).

An incidence structure is nothing more than the idea that two objects coming from distinct classes of things (``points" and ``lines") can be ``incident" with each other. Formally, an \emph{incidence structure} is a triple

\vspace{12pt}
\centerline{($V$, $\mathcal{B}$, $\mathcal{I}$)}
\vspace{12pt}

\noindent where $V$ and $\mathcal{B}$ are any two finite disjoint sets and $\mathcal{I}$ is a symmetric binary relation between $V$ and $\mathcal{B}$, i.e. $\mathcal{I} \subseteq V \times \mathcal{B}$. 

The elements of $V$ are called points, those of $\mathcal{B}$ lines (or blocks), and those of $\mathcal{I}$ incidences (or flags). For convenience, the elements of $V$ are denoted by small-case letters (e.q. $p, q, r, ...$), and those of $\mathcal{B}$ by large-case letters (e.g. $A, B, C, ...$). Using this notation, the elements of $\mathcal{I}$ are pairs, e.g. $(p, A)$, consisting of a point $p$ and a line $A$, indicating that ``p is incident with A", or equivalently, ``p is on A" or ``A passes through p". The influence of geometry is evident in these expressions, although in general, the notions ``point", ``line", and ``incidence" are meant to be broad, and not necessarily restricted to a geometric setting.

\begin{figure}[t!]
\centering
\includegraphics{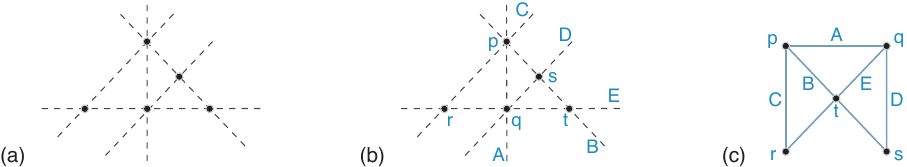}
\caption{}
\label{Figure13}
\end{figure}

Given an arrangement, an incidence structure can be constructed in the following way. Each registration mark and construction line can be labelled with, respectively, small and large letters, in no particular order. The set of incidences can be formed by marking for each point, the lines that pass through that point (equivalently, for each line, the points that are on that line). All the necessary information for this construction (such as, which point is on which line) can be obtained directly from the arrangement itself.

For example, the points and lines of the arrangement in Figure 13(a) can be labelled in the way shown in Figure 13(b). Then, the incidence relation can be obtained by listing every pair of point-line incidence, as in

\vspace{12pt}
\centerline{$\{(p, A), (p, B), (p, C), ..., (r, C), (r, E), ..., (t, B), (t, E)\}$.}
\vspace{12pt}

Incidence structures that represent arrangements for shapes, admit their own graphical presentation in the plane, that is to say, their own drawing. The drawing of an arrangement and the drawing of an incidence structure that represents the arrangement, are not the same thing. Incidence structures cannot capture notions, such as distance, measurement, and betweenness (order of points). As a result, the drawing of an incidence structure need not correspond visually to the drawing of the arrangement, as long as the two drawings represent the same sets of points, lines, and incidences. A drawing for the incidence structure in Figure 13(b), for example, is given in Figure 13(c). Notice that lines in the drawing of an incidence structure are not ``infinitely extending lines". They are line segments of finite extent, and function only as \emph{connections} of points, without further spatial meanings.

An arrangement may also contain construction lines that pass through one point only, such as the arrangement in Figure 14(a). In this case, a construction line is drawn as a line segment, in any orientation or (finite) length; e.g. Figure 14(c). Note that incidence structures cannot capture parallelism. Moreover, potential intersections of lines at places other than the points of the structure, are not accounted for---the elements (points, lines) are drawn using incidence information only.

It is sometimes convenient to use directly the drawing of an arrangement when talking about incidences between construction lines and registration marks, without having to actually draw a corresponding incidence structure separately. In other words, we can informally refer to arrangements as if they were incidence structures, as long as it is clear that we refer to incidences only.

\begin{figure}[t!]
\centering
\includegraphics{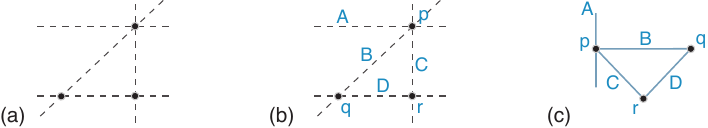}
\caption{}
\label{Figure14}
\end{figure}

\section{Point-line (incidence) geometries and arrangements for shapes}

A finite geometry is a point-line incidence structure that satisfies certain additional rules. The rules that such an incidence structure must satisfy are the minimal ones assumed when speaking about ``a geometry" of points and lines in the plane \cite{Shult2011}. 

A \emph{point-line geometry} is an incidence structure ($V$, $\mathcal{B}$, $\mathcal{I}$), subject to the following rules:

\begin{enumerate}[label=(\roman*)]
    \setlength{\parskip}{0pt}
    \setlength{\itemsep}{0pt plus 1pt}
    \item every line in $\mathcal{B}$ is incident with at least two points of $V$, and 
    \item if two lines in $\mathcal{B}$ are distinct, then there is a point in $V$ incident with one of the lines, which is not incident with the other.
\end{enumerate}

The first rule, attributes to lines the typical characteristic usually associated with them, that is, that a line is formed by ``connecting" at least two points. The second rule, certifies that there are no duplicate lines---distinct lines must be incident with distinct sets of points. The configurations in Figure 5, are two examples of point-line geometries. Point-line geometries are sometimes called spaces. 

From the two rules, it is easy to deduce that any line in a point-line geometry can be identified uniquely with the set of points it is incident with. For this reason, the assumption in the literature is that ``lines" in point-line geometries are subsets of points \cite{BattenCombinatoricsFiniteGeom}.

\begin{figure}[b!]
\centering
\includegraphics[width=\textwidth]{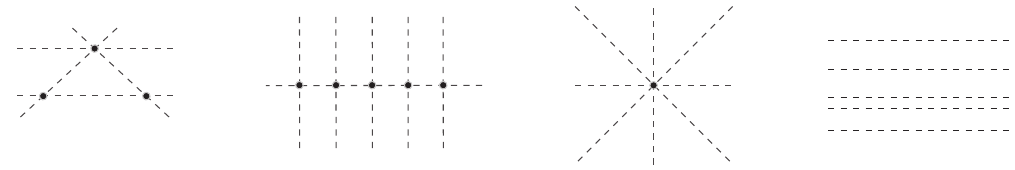}
\caption{Arrangements for shapes that are not point-line geometries.}
\label{Figure15}
\end{figure}

Arrangements for shapes satisfy the second rule automatically. We can never have duplicate construction lines in an arrangement, by definition. However, there are arrangements that do not satisfy the first rule. Arrangements that consist of parallel construction lines or construction lines that pass through one registration mark, do not obey the first rule; see, for example, Figure 15. In either of these two scenarios, the lines cannot be distinguished uniquely on the basis of the points they are incident with---parallel lines (even a single line) have no points to begin with, and having a single point of intersection is not enough for unique identification.

Arrangements, such as those in Figure 15, are not point-line geometries from this strictly mathematical standpoint. One may wish to call them degenerate point-line geometries, although this term will not be used here.

One can further characterize arrangements that do satisfy rules (i) and (ii) in terms of more special types of point-line geometries. 

A point-line geometry is called \emph{near-linear space} (or \emph{partial linear space}) if and only if any two distinct points are incident with at most one line. If we strengthen this rule, so that any two distinct points are incident with \emph{exactly} one line, and further assume that there are at least three non-collinear points, a \emph{linear space} is obtained (the latter rule is a non-triviality condition, since without it, two points connected by a line form a linear space). Linear spaces are thus obtained from near-linear spaces, by adding as lines all pairs of points not already on a line.

Near-linear and linear spaces are well studied in mathematics. They constitute a basic starting point for obtaining other spaces. The classical affine and projective spaces, for example, are linear spaces with extra specialization rules.

\begin{figure}[t!]
\centering
\includegraphics[scale=1.2]{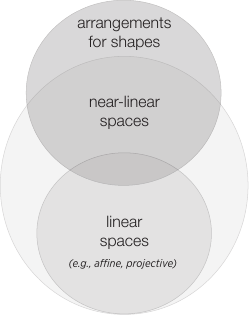}
\caption{}
\label{Figure16}
\end{figure}

The relationship between arrangements for shapes and near-linear and linear spaces can be summarized in a diagram, shown in Figure 16. In particular, there are \emph{some} arrangements for shapes that are near-linear spaces, \emph{some} that are linear spaces, and the rest are arrangements that do not constitute geometries (e.g. Figure 15). For example, all arrangements in Figure 3 are near-linear spaces, the third is also a linear space. In Figure 4, the first three point-line geometries are near-linear spaces, but are not arrangements for shapes (the fourth is not an arrangement, nor a point-line geometry). A longer empirical catalog and an analysis of its properties can be found in the followup study in \cite{HaridisGeomArxiv2020}.

\section{Characterizing determinate and indeterminate shape rules}

Arrangements for shapes and their description as incidence structures provide another way of characterising determinate and indeterminate rules in shape grammars.  

A rule $A \rightarrow B$, where $A$ and $B$ are any two shapes in the universe $U^*$, is \emph{determinate} (i.e. definite transformations can be defined to match the shape $A$) when three conditions are met: (i) three line segments of $A$ do not intersect at a common point, (ii) no two of these line segments are collinear, and (iii) all three are not parallel \cite[p.~261]{StinyShapeBook}. 

Based on (i) through (iii), the arrangement of the shape $A$ must contain at least two registration marks. The smallest arrangement for which this is true is this one

\begin{figure}[h!]
\centering
\includegraphics{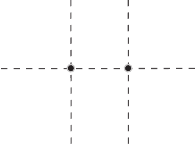}
\label{Figure17}
\end{figure}

\noindent where two construction lines are parallel and a third that intersects both. 

If $A's$ arrangement consists of fewer than two registration marks, then the arrangement is either the empty arrangement, or is of one of the following two types

\vspace{12pt}

\begin{figure}[h!]
\centering
\includegraphics{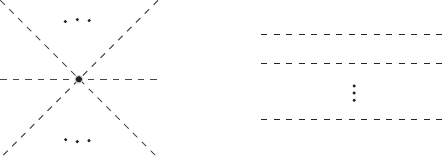}
\label{Figure18}
\end{figure}

\noindent That is, either it has two or more construction lines intersecting at a single registration mark, or it has one or more construction lines parallel to each other. In both of these scenarios, the rule $A \rightarrow B$ is indeterminate (even though there are singular cases where the rule may apply determinately).

A shape $A$ always defines a determinate rule when its underlying arrangement is a point-line geometry. From the two rules for an arrangement in section 2, and the rules for a point-line geometry in section 5, it is easy to see that the smallest possible arrangement for which this happens is the triangle---three registration marks where three construction lines intersect pairwise. There is no arrangement for a shape with fewer than three registration marks that constitutes a point-line geometry.

A shape $A$ may also define a determinate rule when its underlying arrangement is not a point-line geometry. In this case, the arrangement must contain at least two registration marks (for example, the first two arrangements in Figure 15 are not point-line geometries but they define shapes in determinate rules) and the smallest possible arrangement for which this happens is the arrangement shown above with two registration marks and three construction lines (two parallel, one intersecting both).

To decide if a shape $A$ defines a determinate rule, it suffices to construct the arrangement of the shape, then check if the resulting arrangement satisfies the rules for a point-line geometry (it suffices to check if each line is incident with at least two points), in which case the rule is determinate. If not, decide based on the number of registration marks, as explained just above. This construction can be carried out algorithmically in a straightforward way. 

Determinate rules in shape grammars in conjunction with ``description rules" can be used to build rule-based descriptions of design languages. In this case, descriptions can be verbal given in natural language and provide an explanation of the action of a rule when it is applied to generate a new design from a previous one \cite{HaridisSpatialRulesArxiv2021}. The point-line arrangements studied in this paper contribute to the area of mathematics of shapes developed within shape grammar theory. Concurrent studies on the topic of finite topology and topological descriptions for shapes can be found in \cite{HaridisEPB2020} and \cite{HaridisJMA2020}.

\vspace{0.2in}

\section*{Acknowledgement}
\noindent The content of this paper has been revised and subsumed in the following publication:

\vspace{0.2in}

\noindent Alexandros, Haridis. 2024. ``Arrangements containing shapes: mathematical features and their use in visual calculating." In SD Kotsopoulos (Ed.) \emph{Shape Computation: Fifty Years 1972-2022} (Springer Nature). Also in \emph{Mathematics and the Built Environment} book series.



\bibliographystyle{elsarticle-harv}
\bibliography{sample}

\end{document}